# A Fast Algorithm to Calculate Power Sum of Natural Numbers

**Yuyang Zhu**
*Department of Mathematics & Physics, Hefei University, Hefei 230601, P. R. China*
*E-mail: zhuyy@hfuu.edu.cn*

**Abstract:** Permutations can be represented as linear combinations of natural numbers with different powers. In this paper, its coefficient matrix and inverse matrix is derived, and the results show the coefficient matrix is a lower triangular matrix while the inverse matrix is upper triangular. Permutations of n[th] order are used to generate the inverse matrix. The generation function of natural numbers' power sum is derived to calculate the power sum.
**MR(2000) Subject Classification** 11B75, 05A15

## 1. Introduction

The power sum of natural numbers has its applications in number theory and combinatorial mathematics. The existing algorithms of computing power sum includes recursive method, calculus of finite difference, methods of undetermined coefficients and mathematical analysis (see [1-7]). This paper proposes a new algorithm by representing permutation as the linear combinatorial of natural numbers with various powers, and deriving its related coefficient matrix. The matrix is proved to be lower triangular, and its inverse matrix is used to derive the generation function of the power sum. Finally the power sum is computed from the generation function.

## 2. Results

The permutation $P_{n+m}^{n+1} = m(m+1)\cdots(m+n)$ is a $n+1$-th order polynomial of $m$ with integer coefficients, i.e. $P_{n+m}^{n+1} = m^{n+1} + a_1 m^n + \cdots + a_n m$ $(a_1, a_2, \cdots, a_n \in Z)$. By denoting $P_{n+m}^{n+1}$ with $P_{n+1}(m)$, the following result is straightforward to obtain:

**Lemma 2.1** There exists an $n$-th $(n \geq 1)$ order matrix $A_n$ such that
$$(P_1(m), P_2(m), \cdots, P_n(m)) = (m^n, m^{n-1}, \cdots, m)A_n.$$

**Proof** Since polynomials $x, x^2, \cdots, x^n$ are linearly independent in polynomial ring $P$, and $P_1(x), P_2(x), \cdots, P_n(x)$ are 1, 2, $\cdots$, n-th order polynomials of $x$, $P_1(x), P_2(x), \cdots, P_n(x)$ are also linearly independent in polynomial ring $P$. Obviously both are maximally linear independent sets so they can mutually linear expressed, i.e. there exists an $n$-th $(n \geq 1)$ order matrix $A_n$ such that
$$(P_1(x), P_2(x), \cdots, P_n(x)) = (x^n, x^{n-1}, \cdots, x)A_n.$$
Let $x = m$ and the lemma is proved.

**Definition 2.2** Define the $n$-th $(n \geq 1)$ order matrix $A_n$ as $n$-th $(n \geq 1)$ order permutation generation matrix.

**Lemma 2.3** Let $A_n = (a_{ij})_{n \times n}$ be an $n$-th $(n \geq 1)$ order permutation generation matrix, then

(i) If $i + j < n + 1$, then $a_{ij} = 0$;

(ii) If $i + j = n + 1$, then $a_{ij} = 1$;

(iii) If $i + j > n + 1$ and $i < n$, then $a_{ij} = (j-1)a_{i\,j-1} + a_{i+1\,j-1}$;

(iv) If $i = n$, then $a_{n\,i} = (i-1)a_{n\,i-1}$.

**Proof** According to lemma 2.1, $(P_1(m), P_2(m), \cdots, P_n(m)) = (m^n, m^{n-1}, \cdots, m)A_n$, thus

$$\begin{cases} m = P_1(m) = a_{11}m^n + a_{21}m^{n-1} + \cdots + a_{n1}m, \\ m(m+1) = P_2(m) = a_{12}m^n + a_{22}m^{n-1} + \cdots + a_{n2}m, \\ \cdots \qquad \cdots \qquad \cdots \qquad \cdots \\ m(m+1)\cdots(m+n-1) = P_n(m) = a_{1n}m^n + a_{2n}m^{n-1} + \cdots + a_{nn}m. \end{cases} \quad (2.1)$$

According to (2.1), if $i + j < n+1$, then $a_{ij} = 0$, (i) holds.

If $i + j = n+1$, then $a_{ij}$ is an element on the secondary diagonal. According to (2.1),

$$a_{n1} = 1, a_{(n-1)2} = 1, \cdots, a_{1n} = 1,$$

then (ii) holds as well. (iii) and (iv) are proved as follows.

According to (i) and (ii), since $P_j(m) = P_{j-1}(x)(m+j-1)$,

$$P_{j-1}(m) = a_{1(j-1)}m^n + a_{2(j-1)}m^{n-1} + \cdots + a_{(i-1)(j-1)}m^{n-(j-2)} + a_{i(j-1)}m^{n-(j-1)} + \cdots + a_{n(j-1)}m$$
$$= a_{(n+2-j)(j-1)}m^{j-1} + a_{(n+3-j)(j-1)}m^{j-2} + \cdots + a_{n(j-1)}m,$$

Where $a_{(n+2-j)(j-1)} = 1$.

Consider
$$P_j(m) = P_{j-1}(m)(m+j-1) = [a_{(n+2-j)(j-1)}m^{j-1} + a_{(n+3-j)(j-1)}m^{j-2} + \cdots + a_{n(j-1)}m](m+j-1)$$
$$= a_{(n+2-j)(j-1)}m^j + [(j-1)a_{(n+2-j)(j-1)} + a_{(n+3-j)(j-1)}]m^{j-1} + [(j-1)a_{(n+3-j)(j-1)}$$
$$+ a_{(n+4-j)(j-1)}]m^{j-2} + \cdots + [(j-1)a_{(n-1)(j-1)} + a_{n(j-1)}]m^2 + (j-1)a_{n(j-1)}m \quad (2.2)$$

and
$$P_j(m) = a_{(n+1-j)j}m^j + a_{(n+2-j)j}m^{j-1} + \cdots + a_{(n-1)j}m^2 + a_{nj}m. \quad (2.3)$$

By comparing (2.2) and (2.3)
$$a_{(n+1-j)j} = 1, \ a_{(n+2-j)j} = (j-1)a_{(n+2-j)(j-1)} + a_{(n+3-j)(j-1)},$$
$$a_{(n+3-j)j} = (j-1)a_{(n+3-j)(j-1)} + a_{(n+4-j)(j-1)}, \cdots, a_{(n-1)j} = (j-1)a_{(n-1)(j-1)} + a_{n(j-1)},$$
$$a_{nj} = (j-1)a_{n(j-1)}.$$

Thus (iii) and (iv) hold too. The lemma is proved.

It is straightforward to find that $A_2 = \begin{pmatrix} 0 & 1 \\ 1 & 1 \end{pmatrix}$ and according to lemma 2.3, $A_3 = \begin{pmatrix} 0 & 0 & 1 \\ 0 & 1 & 3 \\ 1 & 1 & 2 \end{pmatrix}$. Then

it is trivial to compute the rest:

$$A_4 = \begin{pmatrix} 0 & 0 & 0 & 1 \\ 0 & 0 & 1 & 6 \\ 0 & 1 & 3 & 11 \\ 1 & 1 & 2 & 6 \end{pmatrix}, A_5 = \begin{pmatrix} 0 & 0 & 0 & 0 & 1 \\ 0 & 0 & 0 & 1 & 10 \\ 0 & 0 & 1 & 6 & 35 \\ 0 & 1 & 3 & 11 & 50 \\ 1 & 1 & 2 & 6 & 24 \end{pmatrix}, A_6 = \begin{pmatrix} 0 & 0 & 0 & 0 & 0 & 1 \\ 0 & 0 & 0 & 0 & 1 & 15 \\ 0 & 0 & 0 & 1 & 10 & 85 \\ 0 & 0 & 1 & 6 & 35 & 225 \\ 0 & 1 & 3 & 11 & 50 & 274 \\ 1 & 1 & 2 & 6 & 24 & 120 \end{pmatrix},$$

**Note 1** According to lemma 2.3, $\forall n \in N^+, n \geq 2$, $A_{n-1}$ is a lower left submatrix of $A_n$. The first row of $A_n$ is $a_{11} = a_{12} = \cdots = a_{1(n-1)} = 0, a_{1n} = 1$, and $a_{nn} = (n-1)!$. Thus from known $A_{n-1}$ to compute $A_n$ all we need to compute are $n-2$ elements: $a_{2n}, a_{3n}, \cdots, a_{(n-1)n}$. They are trivial to compute according to (iii) in lemma 2.3.

Similar to the proof of lemma 2.3, we can prove the following theorem:

**Theorem 2.4** Let $A_n^{-1} = (b_{ij})_{n \times n}$ be the inverse of $n$-th order permutation generating matrix, then the following conclusions hold:

(i) if $i+j > n+1$, then $b_{ij} = 0$;

(ii) if $i+j = n+1$, then $b_{ij} = 1$;

(iii) if $i+j \leq n$ and $i \neq 1$, then $b_{ij} = -ib_{i\,j+1} + b_{i-1\,j+1}$;

(iv) if $i=1$, then $b_{1j} = (-1)^{j+n}$.

**Note 2** According to theorem 2.4, $\forall n \in N^+, n \geq 2$, $A_{n-1}^{-1}$ is a upper right submatrix of $A_n^{-1}$. The $n$-th row of $A_n^{-1}$ is
$$b_{2n} = b_{3n} = \cdots = b_{nn} = 0, b_{1n} = 1.$$
Since $b_{11} = (-1)^{1+n}$, to compute $A_n^{-1}$ from a known $A_{n-1}^{-1}$, all we need to find are $n-2$ elements $b_{21}, b_{31}, \cdots, b_{(n-1)1}$ which are straightforward to get from (iii) of theorem 2.4.

**Lemma 2.5** [8] Let $G\{a_n^{(i)}\}$ be the generation function of number sequence $\{a_0^{(i)}, a_1^{(i)}, \cdots, a_n^{(i)}, \cdots\}$ $(i=1,2,\cdots,k)$. Then for any constants $c_1, c_2, \cdots, c_k$,
$$G\{\sum_{i=1}^{k} c_i a_n^{(i)}\} = \sum_{i=1}^{k} c_i G\{a_n^{(i)}\}.$$

**Lemma 2.6** If $(m^n, m^{n-1}, \cdots, m) = (P_1(m), P_2(m), \cdots, P_n(m))A_n^{-1}$, then
$$(G\{m^n\}, G\{m^{n-1}\}, \cdots, G\{m\}) = (G\{P_1(m)\}, G\{P_2(m)\}, \cdots, G\{P_n(m)\})A_n^{-1}.$$

**Proof** Let $A_n^{-1} = \begin{pmatrix} b_{11} & b_{12} & \cdots & b_{1n} \\ b_{21} & b_{22} & \cdots & b_{2n} \\ \cdots & \cdots & \cdots & \cdots \\ b_{n1} & b_{n2} & \cdots & b_{nn} \end{pmatrix}$, and $m^i = b_{1i}P_1(m) + b_{2i}P_2(m) + \cdots + b_{ni}P_n(m)$

$(i=1,2,\cdots,n)$, then according to theorem 2.4, $G\{m^i\} = \sum_{j=1}^{n} b_{ji} G\{P_j(m)\}$, and

$$(G\{m^n\}, G\{m^{n-1}\}, \cdots, G\{m\}) = (\sum_{j=1}^{n} b_{j1}G\{P_j(m)\}, \sum_{j=1}^{n} b_{j2}G\{P_j(m)\}, \cdots, \sum_{j=1}^{n} b_{jn}G\{P_j(m)\})$$

$$= (G\{P_1(m)\}, G\{P_2(m)\}, \cdots, G\{P_n(m)\}) \begin{pmatrix} b_{11} & b_{12} & \cdots & b_{1n} \\ b_{21} & b_{22} & \cdots & b_{2n} \\ \cdots & \cdots & \cdots & \cdots \\ b_{n1} & b_{n2} & \cdots & b_{nn} \end{pmatrix}$$

$$= (G\{P_1(m)\}, G\{P_2(m)\}, \cdots, G\{P_n(m)\})A_n^{-1}.$$

**Theorem 2.7** Let $G\{P_n(m)\}$ be the generation function of $P_{n+m-1}^n = m(m+1)\cdots(m+n-1)$, then
$$G\{P_n(m)\} = \frac{n!x}{(1-x)^{n+1}}.$$

**Proof** In the case of $n=1$,

$$G\{P_1(m)\} = G\{m\} = \sum_{m=1}^{\infty} m x^{m-1} = x \left( \sum_{m=1}^{\infty} x^m \right)' = x \left( \frac{1}{1-x} \right)' = \frac{x}{(1-x)^2}.$$

The theorem holds. Assuming that in the case of $n = k-1$ the theorem also holds, i.e. $G\{P_{k-1}(m)\} = \frac{(k-1)!x}{(1-x)^k}$, then

$$\int_0^x t^{k-2} G\{P_k(m)\} dt = \int_0^x \sum_{m=1}^{\infty} m(m+1)\cdots(m+k-1) t^{m+k-2} dt$$

$$= \sum_{m=1}^{\infty} m(m+1)\cdots(m+k-2) x^{m+k-1} = x^{k-1} G\{P_{k-1}(m)\} = x^{k-1} \frac{(k-1)!x}{(1-x)^k},$$

Take derivatives to both sides of the equation:

$$x^{k-2} G\{P_k(m)\} = \left( x^{k-1} \frac{(k-1)!x}{(1-x)^k} \right)' = \frac{k!x^{k-1}}{(1-x)^{k+1}}.$$

Namely $G\{P_k(m)\} = \frac{k!x}{(1-x)^{k+1}}$, so the theorem also holds in the case of $n = k$. According to the principle of mathematical induction, the theorem is proved.

**Lemma 2.8** [8]  If $b_k = \sum_{i=0}^{k} a_i$, then $G\{b_k\} = \frac{G\{a_k\}}{1-x}$.

### 3. The algorithm to compute power sum of natural numbers

To compute the power sum, the first step is to compute the inverse matrix $A_n^{-1}$ of $n$-th order permutation generating matrix $A_n$, according to theorem 2.4 and mathematical software. The second step is to find the generation function $G\{k^t\}(t = 1, 2, \cdots, n)$ according to lemma 2.6 and theorem 2.7. The third step is to find the generation function of $a_s = \sum_{k=1}^{s} k^t$ according to lemma 2.8. The last step is to calculate the $x^n$ coefficients of generation function $G\{a_s\}$ of $a_s$. The computation is now complete.

**Example 3.1** Compute $\sum_{k=1}^{n} k^4$.

**Solution** Since $A_4 = \begin{pmatrix} 0 & 0 & 0 & 1 \\ 0 & 0 & 1 & 6 \\ 0 & 1 & 3 & 11 \\ 1 & 1 & 2 & 6 \end{pmatrix}$, its inverse $A_4^{-1} = \begin{pmatrix} -1 & 1 & -1 & 1 \\ 7 & -3 & 1 & 0 \\ -6 & 1 & 0 & 0 \\ 1 & 0 & 0 & 0 \end{pmatrix}$, and according to lemma 2.6 and theorem 2.7,

$$G\{k^4\} = (-1)\frac{x}{(1-x)^2} + 7 \cdot \frac{2!x}{(1-x)^3} + (-6)\frac{3!x}{(1-x)^4} + 1 \cdot \frac{4!x}{(1-x)^5} = \frac{(1+11x+11x^2+x^3)x}{(1-x)^5}.$$

Let $a_n = \sum_{k=1}^{n} k^4$, and according to lemma 2.8,

$$G\{a_n\} = \frac{G\{k^4\}}{1-x} = \frac{(1+11x+11x^2+x^3)x}{(1-x)^6} = (x+11x^2+11x^3+x^4)\sum_{k=0}^{\infty} \binom{k+5}{k} x^k.$$

The $x^n$ coefficients are

$$\binom{n+4}{n-1}+11\binom{n+3}{n-2}+11\binom{n+2}{n-3}+\binom{n+1}{n-4}=\frac{n(n+1)(2n+1)(3n^2+3n-1)}{30}.$$

Thus

$$\sum_{k=1}^{n}k^4=\frac{n(n+1)(2n+1)(3n^2+3n-1)}{30}.$$

**Example 3.2** Compute $\sum_{k=1}^{n}k^{12}$.

**Solution** Since $A_{12}^{-1}=$

$$\begin{pmatrix} -1 & 1 & -1 & 1 & -1 & 1 & -1 & 1 & -1 & 1 & -1 & 1 \\ 2047 & -1023 & 551 & -255 & 127 & -63 & 31 & -15 & 7 & -3 & 1 & 0 \\ -86526 & 28510 & -9330 & 3025 & -966 & 301 & -90 & 25 & -6 & 1 & 0 & 0 \\ 611501 & -145750 & 34105 & -7770 & 1701 & -350 & 65 & -10 & 1 & 0 & 0 & 0 \\ -1379400 & 246730 & -42525 & 6951 & -1050 & 140 & -15 & 1 & 0 & 0 & 0 & 0 \\ 1323652 & -179487 & 22827 & -2646 & 266 & -21 & 1 & 0 & 0 & 0 & 0 & 0 \\ -627396 & 63987 & -5880 & 462 & -28 & 1 & 0 & 0 & 0 & 0 & 0 & 0 \\ 159027 & -11880 & 750 & -36 & 1 & 0 & 0 & 0 & 0 & 0 & 0 & 0 \\ -22275 & 1155 & -45 & 1 & 0 & 0 & 0 & 0 & 0 & 0 & 0 & 0 \\ 1705 & -55 & 1 & 0 & 0 & 0 & 0 & 0 & 0 & 0 & 0 & 0 \\ -66 & 1 & 0 & 0 & 0 & 0 & 0 & 0 & 0 & 0 & 0 & 0 \\ 1 & 0 & 0 & 0 & 0 & 0 & 0 & 0 & 0 & 0 & 0 & 0 \end{pmatrix}$$

According to lemma 2.6 and theorem 2.7,

$$G\{k^{12}\}=(-1)\frac{x}{(1-x)^2}+2047\frac{2!x}{(1-x)^3}+(-86526)\frac{3!x}{(1-x)^4}+(611501)\frac{4!x}{(1-x)^5}+(-1379400)\frac{5!x}{(1-x)^6}$$

$$+(1323652)\frac{6!x}{(1-x)^7}+(-627396)\frac{7!x}{(1-x)^8}+159027\frac{8!x}{(1-x)^9}+(-22275)\frac{9!x}{(1-x)^{10}}$$

$$+1705\frac{10!x}{(1-x)^{11}}+(-66)\frac{11!x}{(1-x)^{12}}+\frac{12!x}{(1-x)^{13}}$$

$$=\frac{1}{(1-x)^{13}}[x(1+4083x+478271x^2+10187685x^3+66318474x^4+162512286x^5$$

$$+162512286x^6+66318474x^7+10187685x^8+478271x^9+4083x^{10}+x^{11})].$$

Let $a_n=\sum_{k=1}^{n}k^{12}$, and according to lemma 2.8,

$$G\{a_n\}=\frac{G\{k^{12}\}}{1-x}=\frac{1}{(1-x)^{14}}[x(1+4083x+478271x^2+10187685x^3+66318474x^4+162512286x^5$$

$$+162512286x^6+66318474x^7+10187685x^8+478271x^9+4083x^{10}+x^{11})]$$

$$=(x+4083x^2+478271x^3+10187685x^4+66318474x^5$$

$$+162512286x^6+162512286x^7+66318474x^8$$

$$+10187685x^9+478271x^{10}+4083x^{11}+x^{12})\sum_{k=0}^{\infty}\binom{k+13}{k}x^k.$$

The $x^n$ coefficients on the right hand side are

$$\binom{n+12}{n-1}+4083\binom{n+11}{n-2}+478271\binom{n+10}{n-3}+10187685\binom{n+9}{n-4}+$$
$$66318474\binom{n+8}{n-5}+162512286\binom{n+7}{n-6}+162512286\binom{n+6}{n-7}+66318474\binom{n+5}{n-8}$$
$$+10187685\binom{n+4}{n-9}+478271\binom{n+3}{n-10}+4083\binom{n+2}{n-11}+\binom{n+1}{n-12}$$
$$=-\frac{691n}{2730}+\frac{5n^3}{3}-\frac{33n^5}{10}+\frac{22n^7}{7}-\frac{11n^9}{6}+n^{11}+\frac{n^{12}}{2}+\frac{n^{13}}{13}.$$

Thus $\sum_{k=1}^{n}k^{12}=-\frac{691n}{2730}+\frac{5n^3}{3}-\frac{33n^5}{10}+\frac{22n^7}{7}-\frac{11n^9}{6}+n^{11}+\frac{n^{12}}{2}+\frac{n^{13}}{13}$.

**Note 3** According to Note 2, If the $n$-th order inverse of permutation generation matrix $A_n^{-1}$, then for arbitrary natural number $m$ ($1\leq m\leq n$), $A_m^{-1}$ is upper right triangular submatrix of $A_n^{-1}$. Thus if $A_n^{-1}$ has been computed, then the power sum $\sum_{i=1}^{k}i^m$ can be computed for arbitrary $m$ ($1\leq m\leq n$). For example, $\sum_{i=1}^{k}i^{12}$ can be calculated by the first column of $A_{12}^{-1}$, $\sum_{i=1}^{k}i^{11}$ can be calculated by the second column of $A_{12}^{-1}$, $\sum_{i=1}^{k}i^{10}$ can be calculated by the third column of $A_{12}^{-1}$, $\sum_{i=1}^{k}i^{9}$ can be calculated by the fourth column of $A_{12}^{-1}$, and $\sum_{i=1}^{k}i^{4}$ can be calculated by the ninth column of $A_{12}^{-1}$. This algorithm is straightforward to implement on the computer, i.e. $\sum_{i=1}^{k}i^{50}$ can be calculated in a few minutes on Mathematics.